\documentclass[11pt]{article}
\usepackage{graphicx}
\usepackage{amsmath}
\usepackage[dvips]{epsfig}
\usepackage{amssymb}
\begin{document}
\begin{center} {\Large \bf  $q$-Bernstein
polynomials,  $q$-Stirling numbers and  $q$-Bernoulli polynomials}
\\ \vspace*{12 true pt}  T.  Kim
\vspace*{12 true pt} \\
 Division of General Education-Mathematics,\\
 Kwangwoon University, Seoul 139-701,  Korea \\
                        \end{center}
%--------------------------------------------------------------------
\vspace*{12 true pt} \noindent {\bf Abstract :} In this paper, we
give new identities involving Phillips $q$-Bernstein polynomials
and we derive some interesting properties of $q$-Bernstein
polynomials associated with $q$-Stirling numbers and $q$-Bernoulli
polynomials.

\vspace*{12 true pt} \noindent {\bf 2000 Mathematics Subject
Classification :} 11B68, 11S40, 11S80

 \vspace*{12 true pt} \noindent {\bf
Key words :} $q$-Bernoulli polynomials, $q$-Bernstein polynomials,
$q$-Stirling numbers

\begin{center} {\bf 1. Introduction } \end{center}

When one talks of $q$-extension, $q$ is variously considered  as
an indeterminate, a complex number $q\in \mathbb{C},$ or $p$-adic
number $q\in\Bbb C_p .$ If $q\in \Bbb C$, then we always assume
that $|q|<1.$  If $q\in  \mathbb{C}_p,$ we usually assume that
$|1-q|_p< 1$. Here, the symbol $| \cdot|_p $ stands for the
$p$-adic absolute value on  $\mathbb{C}_p$ with $|p|_p \leq
{1}/{p}.$ For each $x$, the $q$-basic numbers are defined by
$$[x]_q =\frac{1-q^x}{1-q}, \mbox{ and } [n]_q!=[n]_q [n-1]_q \cdots [2]_q [1]_q, n \in \Bbb N, \text{ (see [1-17])}.$$
Throughout this paper we assume that $ q \in \Bbb C$ with $|q|<1$
and we use the notation of Gaussian binomial coefficient in the
form
$${\binom{n}{k}}_q= \dfrac{[n]_q !}{[k]_q![n-k]_q!}=\dfrac{[n]_q [n-1]_q \cdots [n-k+1]_q}{[k]_q!},  n, k \in \Bbb N.$$
Note that $$ \lim_{q \rightarrow
1}\binom{n}{k}_q=\binom{n}{k}=\dfrac{n(n-1)\cdots (n-k+1)}{k!},
\text { (see [4-12])} .$$
 The Gaussian binomial coefficient
satisfies the following recursion formula:
$$\binom{n+1}{k}_q=\binom{n}{k-1}_q+ q^k \binom{n}{k}_q= q^{n-k} \binom{n}{k-1}_q+\binom{n}{k}_q, \text{ (see [7, 8])}.\eqno(1)$$
The $q$-binomial formulae are known as
$$(1-b)_q^n=(b:q)_n=\prod_{i=1}^n (1-bq^{i-1})= \sum_{i=0}^n {\binom{n}{i}}_q q^{{\binom{i}{2}}} (-1)^ib^i, \eqno(2) $$
and
$$\dfrac{1}{(1-b)_q^n}=\dfrac{1}{(b:q)_n}=\dfrac{1}{\prod_{i=1}^n (1-bq^{i-1})}
= \sum_{i=0}^\infty  {\binom{n+i-1}{i}}_q   b^i, \text{ (see [7,
8])}. $$

Now, we define the $q$-exponential function as follows:
$$ \lim_{n \rightarrow
\infty } \dfrac{1}{(x:q)_n}= \lim_{n \rightarrow \infty }
\sum_{k=0}^\infty  \binom{n+k-1}{k}_q x^k =\sum_{k=0}^\infty
\dfrac{ x^k (1-q)^k}{[k]_q!}=e_q(x(1-q)). \eqno(3)$$
A Bernoulli trial involves performing an experiment once and noting whether
 a particular event $A$ occurs. The outcome of Bernoulli trial is said to be ``success"  if $A$ occurs and a ``failure" otherwise.
Let $k$ be the number of successes in $n$ independent Bernoulli trials, the probabilities of $k$ are given by the binomial probability law:
$$p_n(k)=\binom{n}{k}p^k(1-p)^{n-k},  \text{ for $k=0, 1, \cdots, n$,}$$
where $p_n(k)$ is the probability of $k$ successes in $n$ trials. For example, a communication system transmit binary information over channel
that introduces random bit errors with probability $\xi=10^{-3}$. The transmitter transmits each information bit three times, an a decoder
takes a majority vote of the received bits to decide on what the transmitted bit was. The receiver can correct a single error, but it will make
the wrong decision if the channel introduces two or more errors. If we view each transmission as a Bernoulli trial in which a ``success" corresponds to the introduction of an error, then the probability of two or more errors in three Bernoulli trial is
$$p (k\geq 2)=\binom{3}{2}(0.001)^2(0.999)+\binom{3}{3}(0.001)^3\approx 3(10^{-6}), \text{ see [18]}.$$
Let $C[0,1]$ denote the set of continuous function on $[0, 1]$.
For $f\in C[0, 1]$, Bernstein introduced the following well known linear operator in [2]:
$$B_n(f|x)=\sum_{k=0}^nf(\frac{k}{n})\binom{n}{k}x^k(1-x)^{n-k}=\sum_{k=0}^nf(\frac{k}{n})B_{k,n}(x).$$
Here $B_n(f|x)$ is called the Bernstein operator of order $n$ for $f$. For $k, n \in \Bbb Z_{+}$, the Bernstein polynomials of degree $n$ is defined 
by $$B_{k,n}(x)=\binom{n}{k}x^k (1-x)^{n-k}.$$ By the definition of Bernstein polynomials, we can see that
Bernstein basis is the probability mass function of binomial distribution. 
Based on the $q$-integers Phillips introduced the $q$-analogue of
well known Bernstein polynomials (see [15, 16]).
For $ f \in C[0,1],$ Phillips introduced the $q$-extension of $\Bbb B_n( f|x)$ as follows:
$$ \Bbb B_{n, q}(f \mid x)=\sum_{k=0}^n B_{k,n}(x, q) f \left( \dfrac{[k]_q}{[n]_q}\right)
=\sum_{k=0}^n  f \left( \dfrac{[k]_q}{[n]_q}\right) \binom{n}{k}_q
x^k (1-x)_q^{n-k}. \eqno(4)$$ Here $ \Bbb B_{n, q}(f \mid x)$ is
called the $q$-Bernstein operator of order $n$ for $f$.
For $k, n \in \Bbb Z_+,$ the $q$-Bernstein polynomial of degree
$n$ is defined by
$$B_{k,n}(x,q)=  {\binom{n}{k}}_q x^k (1-x)_q^{n-k}, x \in [0, 1]. \eqno(5) $$
For example, $B_{0,1}(x,q)=1-x,  B_{1,1}(x,q)=x, $ and
$B_{0,2}(x,q)=1-[2]_qx +qx^2, \cdots. $ Also $B_{k,n}(x,q)=0$ for
$k>n,$ because $\binom{n}{k}_q=0$.
 The $q$-binomial distribution associated with the $q$-boson oscillator are introduced in [19, 20]. For $n, k \in \Bbb Z_{+},$ its probabilities are given by
$$p(X=k)=\binom{n}{k}_q x^k(1-x)_q^{n-k}, \text{ where $x\in [0, 1]$}.$$  This distributions are  studied  by several authors and has applications in physics as well as  in approximation theory due to the $q$-Bernstein polynomials and the $q$-Bernstein operators (see [16, 19, 20]).
From the definition of $q$-Bernstein polynomials, we note that the $q$-Bernstein basis is the probability mass function 
of $q$-binomial distribution.
 Recently, several authors have
studied the analogs of Bernstein polynomials (see [1, 2, 5, 8, 9,
10, 15, 16, 17]). In [5], Gupta-Kim-Choi-Kim gave the generating
function of Phillips $q$-Bernstein polynomials as follows:
$$ \aligned   \dfrac{x^k t^k}{[k]_q!} e_q((1-x)_q t) &=\dfrac{x^k t^k}{[k]_q!} \sum_{n=0}^\infty \dfrac{ (1-x)_q^n
t^n}{[n]_q!}\\
&  = \sum_{n=k}^\infty \binom nk_q   \dfrac{ x^k
(1-x)_q^{n-k}}{[n]_q!}t^n  \\
& = \sum_{n=k}^\infty B_{k,n}(x, q) \dfrac{ t^n}{[n]_q!}.
\endaligned
$$
Because $B_{k,0}(x,q)=B_{k,1}(x,q)=B_{k,2}(x,q) =\cdots
=B_{k,k-1}(x,q)=0,$ we obtain the generating function for
$B_{k,n}(x,q)$ as follows:
$$F_q^{(k)}(t, x)=\dfrac{x^k t^k}{[k]_q!} e_q((1-x)_q t)=\sum_{n=0}^\infty B_{k,n}(x, q) \dfrac{ t^n}{[n]_q!}, \text{ see [5]},$$
where $n,  k \in \Bbb Z_+$ and $ x \in [0, 1]$.

Notice that
$$B_{k, n}(x, q) =\left \{\begin{array}{ll}
 {\binom{n}{k}}_q x^k (1-x)_q^{n-k},   & \mbox{ if } n \geq k \\
0, &  \mbox{ if } n < k ,
\end{array} \right.
$$
for $n,  k \in \Bbb Z_+$ (see [5, 15, 16]).

In this paper we study the  generating function of Phillips
$q$-Bernstein polynomial and give some identities on the Phillips
$q$-Bernstein polynomials. From the generating function of
$q$-Bernstein polynomial, we derive recurrence relation and
derivative of the Phillips $q$-Bernstein polynomials. Finally, we
investigate some interesting properties of  $q$-Bernstein
polynomials related to $q$-Stirling numbers and  $q$-Bernoulli
polynomials.

\bigskip

\begin{center} {\bf 2.  $q$-Bernstein
polynomials,  $q$-Stirling numbers and  $q$-Bernoulli polynomials}
\end{center}

Let
$$F_q^{(k)}(t, x)=\sum_{n=0}^\infty B_{k,n}(x, q) \dfrac{
t^n}{[n]_q!}.$$ From (5) and (3), we note that

$$ \aligned  F_q^{(k)}(t, x) & = \sum_{n=0}^\infty \binom nk_q  x^k
(1-x)_q^{n-k} \dfrac{ t^n}{[n]_q!}  \\
&=\sum_{n=0}^\infty \binom {n+k}{k}_q  \dfrac{  x^k (1-x)_q^{n} }{[n+k]_q!} t^{n+k} \\
&=\dfrac{x^k t^k}{[k]_q!} \sum_{n=0}^\infty \dfrac{ (1-x)_q^{n}
}{[n]_q!} t^{n} \\
& = \dfrac{x^k t^k}{[k]_q!} e_q((1-x)_q t),
\endaligned
$$
where $ n,  k \in \Bbb Z_+$ and $ x \in [0, 1].$

Note that
$$ \lim_{ q \rightarrow 1}F_q^{(k)}(t, x)= \dfrac{x^k t^k}{k!}e^{(1-x)t}= \sum_{n=0}^\infty  B_{k, n}(x) \dfrac{t^n}{n!},$$
where $ B_{k, n}(x) $ are the Bernstein polynomial of degree $n$.

The $q$-derivative $D_qf$ of function $f$ is defined by
$$(D_q f)(x)=\dfrac{df(x)}{d_qx}= \dfrac{f(x)-f(qx)}{ (1-q)x}, \text{ (see [6])}. \eqno(7)$$
From (7), we note that
$$D_q(fg)(x)=g(x) D_qf(x)+ f(qx)D_q g(x), \text{ (see [6])}. \eqno(8)$$
The  $q$-Bernstein operator is given by
$$ \Bbb B_{n, q}(f \mid x)=\sum_{k=0}^n B_{k,n}(x, q) f \left( \dfrac{[k]_q}{[n]_q}\right), \text{ (see Eq. (4))}.$$
Thus, we have
$$ \Bbb B_{n, q}(1\mid x)=\sum_{k=0}^n B_{k,n}(x, q) =\sum_{k=0}^n {\binom{n}{k}}_q x^k (1-x)_q^{n-k}=1,$$
and
$$ \Bbb B_{n, q}(x \mid x)=\sum_{k=0}^n \left( \dfrac{[k]_q}{[n]_q}\right)
B_{k,n}(x, q) =x \sum_{k=0}^{n-1} {\binom{n-1}{k}}_q x^k
(1-x)_q^{n-k}=x,$$ where  $ x \in [0, 1]$ and $ n,  k \in \Bbb
Z_+$.

For $f \in C[0,1]$, we have

$$ \aligned \Bbb B_{n, q}(f \mid x) & = \sum_{k=0}^n  f \left( \dfrac{[k]_q}{[n]_q}\right)
B_{k,n}(x, q) \\
&= \sum_{k=0}^n f \left( \dfrac{[k]_q}{[n]_q}\right)  {\binom{n}{k}}_q x^k (1-x)_q^{n-k} \\
&= \sum_{k=0}^n f \left( \dfrac{[k]_q}{[n]_q}\right) x^k
{\binom{n}{k}}_q  \sum_{j=0}^{n-k} \binom{n-k}{j}_q (-1)^j
q^{{\binom{j}{2}}}  x^j .
\endaligned
$$
It is easy to show that
$$\binom{n}{k}_q \binom{n-k}{j}_q =\binom{n}{k+j}_q \binom{k+j}{k}_q. $$
Let $k+j=m.$ Then we have
$$\binom{n}{k}_q \binom{n-k}{j}_q =\binom{n}{m}_q \binom{m}{k}_q. \eqno(10) $$
By (9) and (10), we easily get

$$  \Bbb B_{n, q}(f \mid x) = \sum_{m=0}^n {\binom{n}{m}}_q x^m
  \sum_{k=0}^{m} \binom{m}{k}_q q^{{\binom{m-k}{2}}} (-1)^{m-k}
 f \left( \dfrac{[k]_q}{[n]_q}\right). \eqno(11)
$$
Therefore, we obtain the following proposition.

\bigskip
{ \bf Proposition 1.} For $f \in C[0,1]$ and $ n \in
\mathbb{Z}_+,$ we have
$$  \Bbb B_{n, q}(f \mid x) = \sum_{m=0}^n {\binom{n}{m}}_q x^m
  \sum_{k=0}^{m} \binom{m}{k}_q q^{{\binom{m-k}{2}}} (-1)^{m-k}
 f \left( \dfrac{[k]_q}{[n]_q}\right). \eqno(11)
$$
\bigskip

It is well known that the second kind Stirling numbers are defined
by
$$\dfrac{(e^t-1)^k}{k!}=\dfrac{1}{k!} \sum_{l=0}^k \binom kl (-1)^{k-l} e^{lt} = \sum_{n=0}^\infty  S(n, k) \dfrac{t^n}{n!},
\eqno(12)$$ where $ k \in \Bbb N$ (see [7, 8, 9, 10, 17]).

Let $\Delta$ be the shift difference operator with $ \Delta
f(x)=f(x+1)-f(x)$. By iterative process, we see that

$$\Delta^n f(0)= \sum_{k=0}^n \binom nk (-1)^{n-k}f(k), \text{ for } n \in \Bbb N. \eqno(13) $$

From (12) and (13), we have
$$ \aligned \sum_{n=0}^\infty  S(n, k) \dfrac{t^n}{n!} & =\dfrac{1}{k!} \sum_{l=0}^k \binom kl (-1)^{k-l} e^{lt} \\
&= \sum_{n=0}^\infty \left(\dfrac{1}{k!} \sum_{l=0}^k \binom kl
(-1)^{k-l}l^n \right) \dfrac{t^n}{n!}\\
&= \sum_{n=0}^\infty  \dfrac{\Delta^k 0^n }{k!} \dfrac{t^n}{n!},
\text{ (see [7, 8, 9])}.
\endaligned
\eqno(14)
$$
By comparing the coefficients on the both sides of (14), we get
$$S(n, k)= \dfrac{\Delta^k 0^n }{k!} , \text{ for }n, k \in \Bbb Z_+. \eqno(15)$$
Now, we consider the $q$-extension of (13).  Let $(Eh)(x)=h(x+1)$
be the shift operator. Then the $q$-difference operator is defined
by
$$\Delta_q^n := (E-I)_q^n= \prod_{i=1}^n (E-Iq^{i-1}), \text{ (see [7])}, $$
where $I$ is an identity operator.

For $f \in C[0,1]$ and $n \in \Bbb N$, we have

$$\Delta_q^n f(0)=  \sum_{k=0}^{n} \binom{n}{k}_q   (-1)^{k} q^{{\binom{k}{2}}}  f (n-k )=
 \sum_{k=0}^{n} \binom{n}{k}_q   (-1)^{n-k} q^{{\binom{n-k}{2}}}  f (k ). \eqno(16)$$
Note that (16) is exactly $q$-extension of (13). That is, $
\lim_{q \rightarrow 1}\Delta_q^n f(0)=\Delta^n f(0).$

 As the
$q$-extension of (12), the second kind $q$-Stirling numbers are
defined by
$$\dfrac{q^{-{\binom{k}{2}}}}{[k]_q!} \sum_{j=0}^{k} (-1)^{k-j} \binom{k}{j}_q  q^{{\binom{k-j}{2}}} e^{[j]_qt}
= \sum_{n=0}^\infty  S(n, k: q) \dfrac{t^n}{n!}, \text{ (see [7,
8])}. \eqno(17).
$$
By (16), we obtain the following theorem.

\bigskip
{ \bf Theorem 2.} For $f \in C[0,1]$ and $ n \in \mathbb{Z}_+,$ we
have
$$  \Bbb B_{n, q}(f \mid x) = \sum_{k=0}^n {\binom{n}{k}}_q x^k
   \Delta_q^k f \left( \dfrac{0}{[n]_q}\right).
$$
 \bigskip

In the special case, $f(x)=x^m ( m \in \Bbb Z_+),$ we obtain the
following corollary.

\bigskip
{ \bf Corollary 3.} For $ x \in [0,1]$ and $ m, n \in
\mathbb{Z}_+,$ we have
$$  [n]_q^m  \Bbb B_{n, q}(x^m \mid x) = \sum_{k=0}^n {\binom{n}{k}}_q x^k
   \Delta_q^k 0^m.
$$
 \bigskip

By (17),  we easily get
$$ \aligned  S(n, k : q) &
 =\dfrac{q^{-{\binom{k}{2}}}}{[k]_q!} \sum_{j=0}^{k} (-1)^{j} q^{{\binom{j}{2}}} \binom{k}{j}_q [k-j]_q^n\\
&= \dfrac{q^{-{\binom{k}{2}}}}{[k]_q!} \sum_{j=0}^{k} (-1)^{k-j} q^{{\binom{k-j}{2}}} \binom{k}{j}_q [j]_q^n\\
&= \dfrac{q^{-{\binom{k}{2}}}}{[k]_q!}\Delta_q^k 0^m.
\endaligned
\eqno(18)
$$

The equation (18) seems to be the $q$-extension of the equation
(15). That is,  $ \lim_{q \rightarrow 1}S(n, k : q)=S(n, k ).$

By Corollary 3 and (18), we obtain the following corollary.

\bigskip
{ \bf Corollary 4.} For $ x \in [0,1]$ and $ m, n \in
\mathbb{Z}_+,$ we have
$$  [n]_q^m  \Bbb B_{n, q}(x^m \mid x) = \sum_{k=0}^n {\binom{n}{k}}_q x^k
  [k]_q! q^{{\binom{k}{2}}} S(m, k :q) .
$$
 \bigskip
From (1) and (5), for $ 0 \leq k \leq n,$ we have
$$ \aligned  & q^k (1- q^{n-k-1} x) B_{k, n-1}(x,q)+ x B_{k-1, n-1}(x,q)\\
&= q^k (1- q^{n-k-1} x)\binom{n-1}{k}_q x^k (1-x)_q ^{n-1-k} + x  \binom{n-1}{k-1}_q x^{k-1} (1-x)_q^{n-k}\\
&= q^k \binom{n-1}{k}_q x^k (1-x)_q ^{n-k} +  \binom{n-1}{k-1}_q
x^{k} (1-x)_q^{n-k}\\
&= \binom{n}{k}_q x^{k} (1-x)_q^{n-k} .
\endaligned
\eqno(19)
$$
By (2), (7) and (8), we get
$$ \dfrac{d B_{k, n}(x, q)}{d_qx}= -\binom nk_q x^k [n-k]_q ( 1-qx)_q^{n-k-1}+
\binom nk_q [k]_q x^{k-1} ( 1-qx)_q^{n-k}. \eqno(20) $$

From the definition of Gaussian binomial coefficient (=
$q$-binomials coefficient) and (2), we note that
$$\binom nk_q [k]_q x^{k-1} ( 1-qx)_q^{n-k} = q^{-(k-1)} [n]_q B_{k-1, n-1}(qx, q), \eqno(21) $$
and
$$\binom nk_q x^{k} [n-k]_q  ( 1-qx)_q^{n-k-1} = [n]_q q^{-k} B_{k, n-1}(qx, q). $$

By (20) and (21), we see that
$$ \dfrac{d B_{k, n}(x, q)}{d_qx}= [n]_q q^{-k} ( q B_{k-1, n-1}(qx, q)- B_{k, n-1}(qx, q)). \eqno(22) $$
\bigskip

Thus, we note that the $q$-derivative of the $q$- Bernstein
polynomials of degree $n$ are also  polynomial of degree $n-1$.
Therefore, by (19) and (22), we obtain the following recurrence
formulae:

\bigskip
{ \bf Theorem 5}(Recurrence formulae for $B_{k, n}(x, q)$).
 For
$k, n \in \mathbb{Z}_+$  and $ x \in [0, 1]$,  we have
$$  q^k (1- q^{n-k-1} x) B_{k, n-1}(x,q)+ x B_{k-1, n-1}(x,q)= B_{k, n}(x,q),$$
and
$$ \dfrac{d B_{k, n}(x, q)}{d_qx}= [n]_q q^{-k} ( q B_{k-1, n-1}(qx, q)- B_{k, n-1}(qx, q)). $$
 \bigskip

We also get from (5) and (6) that

$$ \aligned
&\dfrac{[n-k]_q}{[n]_q} B_{k, n}(x, q)+  \dfrac{[k+1]_q}{[n]_q}
B_{k+1, n}(x, q) \\
& =(1- x q^{n-k-1} )\binom{n-1}{k}_q x^k (1-x)_q ^{n-k-1} + x  \binom{n-1}{k}_q x^{k} (1-x)_q^{n-k-1}\\
&= (1- x q^{n-k-1} )B_{k, n-1}(x, q) + x  B_{k, n-1}(x, q)\\
&= B_{k, n-1}(x, q)+ x[n-k-1]_q (1-q) B_{k, n-1}(x, q).
\endaligned
\eqno(23)
$$
By (23), we obtain the following theorem.

\bigskip
{ \bf Theorem 6.}  For $k, n \in \mathbb{Z}_+$  and $ x \in [0,
1]$,  we have
$$   \dfrac{[n-k]_q}{[n]_q} B_{k, n}(x, q)+  \dfrac{[k+1]_q}{[n]_q}
B_{k+1, n}(x, q) =B_{k, n-1}(x, q)+ x[n-k-1]_q (1-q) B_{k, n-1}(x,
q).
$$
 \bigskip

From Theorem 6 we note that $q$-Bernstein polynomials can be
written as a linear combination of polynomials of higher order.

For $k, n \in \Bbb N, $ we easily get from (5) that $q$-Bernstein
polynomials can be expressed in the form
$$ \aligned
&\dfrac{[n-k+1]_q}{[k]_q}  \left(  \dfrac{x}{1-xq^{n-k} } \right)
x^{k-1} (1-x)_q^{n-k+1} \binom {n}{k-1}_q  \\
& = \dfrac{[n]_q !}{[k]_q ! [n-k]_q!}  x^k (1- x)_q^{n-k}\\
&= \binom{n}{k}_q x^k (1-x)_q ^{n-k}\\
&=  B_{k, n}(x, q).
\endaligned
\eqno(24)
$$
By (24), we obtain the following proposition.

\bigskip
{ \bf Proposition 7.}  For $n, k \in \mathbb{N}$  and $ x \in [0,
1]$, we have
$$    B_{k, n}(x, q) =\dfrac{[n-k+1]_q}{[k]_q}  \left(  \dfrac{x}{1-xq^{n-k} } \right) B_{k-1, n}(x,
q).
$$
 \bigskip

The $q$-Bernstein polynomials of degree $n$ can be written in
terms of power basis $\{ 1, x, x^2, \cdots, x^n \}.$  By using the
definition of $q$-Bernstein polynomial and $q$-binomial theorem,
we get

$$ \aligned  B_{k, n}(x, q) & =  {\binom{n}{k}}_q  x^k (1-x)_q^{n-k}
 ={\binom{n}{k}}_q  x^k \sum_{i=0}^{n-k} {\binom{n-k}{i}}_q (-1)^i q^{ \binom{i}{2}} x^i  \\
&= \sum_{i=0}^{n-k} {\binom{n-k}{i}}_q  {\binom{n}{k}}_q (-1)^i q^{ \binom{i}{2}} x^{i+k}\\
&= \sum_{i=k}^{n} {\binom{n-k}{i-k}}_q  {\binom{n}{k}}_q
(-1)^{i-k}  q^{ \binom{i-k}{2}} x^{i}.
\endaligned
\eqno(25)
$$
By simple calculation, we easily see that
$$\binom{n}{k}_q \binom{n-k}{i-k}_q =\binom{n}{i}_q \binom{i}{k}_q. \eqno(26) $$
Therefore, by (25) and (26), we obtain the following theorem.

\bigskip
{ \bf Theorem 8.}  For $k, n \in \mathbb{Z}_+$  and $ x \in [0,
1]$,  we have
$$   B_{k, n}(x, q)=\sum_{i=k}^{n} {\binom{n}{i}}_q  {\binom{i}{k}}_q
(-1)^{i-k}  q^{ \binom{i-k}{2}} x^{i}.
$$
 \bigskip

We get from the properties of $q$-Bernstein polynomials that

$$ \aligned  \sum_{k=1}^{n} \dfrac{{\binom{k}{1}}_q }{{\binom{n}{1}}_q } B_{k, n}(x, q) & =
\sum_{k=1}^{n} \dfrac{[k]_q }{[n]_q} {\binom{n}{k}}_q x^k
(1-x)_q^{n-k} \\
& =\sum_{k=1}^{n}  {\binom{n-1}{k-1}}_q x^k
(1-x)_q^{n-k} \\
&=  x \sum_{k=0}^{n-1} {\binom{n-1}{k}}_q  x^k (1-x)_q^{n-k-1}=x,
\endaligned
$$
and that
$$ \aligned  \sum_{k=2}^{n} \dfrac{{\binom{k}{2}}_q }{{\binom{n}{2}}_q } B_{k, n}(x, q) &
=\sum_{k=2}^{n}  {\binom{n-2}{k-2}}_q x^k
(1-x)_q^{n-k} \\
&=  x^2 \sum_{k=0}^{n-2} {\binom{n-2}{k}}_q  x^k
(1-x)_q^{n-k-2}=x^2.
\endaligned
$$
Continuing this process, we obtain
$$\sum_{k=i}^{n} \dfrac{{\binom{k}{i}}_q }{{\binom{n}{i}}_q } B_{k, n}(x, q)=x^i.$$
Therefore, we obtain the following theorem.

\bigskip
{ \bf Theorem 9.}  For $k, i \in \mathbb{Z}_+$  and $ x \in [0,
1]$,  we have
$$\sum_{k=i}^{n} \dfrac{{\binom{k}{i}}_q }{{\binom{n}{i}}_q } B_{k, n}(x, q)=x^i.$$
 \bigskip

Now we define $q$-Bernoulli polynomials of order $k$ as follows:

$$ \left(  \dfrac{z}{e^z-1}  \right)^k e_q(zx)= \sum_{n=0}^\infty \beta_n^{(k)}(x, q) \dfrac{z^n}{[n]_q!}, \quad k \in \Bbb
N. \eqno(27)$$ From the generating function (27) of $q$-Bernoulli
polynomials and (3), we derive

$$ \aligned \left(  \dfrac{z}{e^z-1}  \right)^k e_q(zx) & =
\left(  \sum_{m=0}^{\infty} B_m^{(k)} \dfrac{z^m}{m!}\right )  \left(  \sum_{l=0}^{\infty} \dfrac{x^lz^l}{[l]_q!}\right)\\
&=  \sum_{n=0}^{\infty} \left(  \sum_{m=0}^{n} \dfrac{B_m^{(k)}
x^{n-m} [n]_q!}{m! [n-m]_q!} \right ) \dfrac{z^n}{[n]_q!} \\
&=\sum_{n=0}^{\infty} \left(  \sum_{m=0}^{n} \dfrac{ [m]_q!}{m!}
B_m^{(k)} \binom nm_q x^{n-m} \right ) \dfrac{z^n}{[n]_q!} ,
\endaligned \eqno(28)
$$
where $B_m^{(k)}$ are the $n$-th  Bernoulli numbers of order
$k$(see [6]).

From (27) and (28), we easily get
$$\beta_n^{(k)}(x, q)=  \sum_{m=0}^{n}  \binom nm_q \dfrac{ [m]_q!}{m!}
 x^{n-m} B_m^{(k)}, \eqno(29)$$
 where $B_m^{(k)}$ are  the $m$-th  ordinary  Bernoulli numbers of order
$k$.

From (26) and (27), we note that

$$ \aligned  \dfrac{(tx)^k}{ [k]_q!} e_q((1-x)_q t)
&= \dfrac{x^k(e^t-1)^k }{ [k]_q!} \left(  \dfrac{t}{e^t-1}  \right)^k e_q((1-x)_q t) \\
& =\dfrac{k! }{ [k]_q!} x^k  \left(  \sum_{m=0}^{\infty} S(m, k)
\dfrac{t^m}{m!}  \right ) \left( \sum_{n=0}^{\infty}
\beta_n^{(k)}((1-x)_q, q)\dfrac{t^n}{[n]_q!}  \right )\\
&= \dfrac{k! }{ [k]_q!} x^k \sum_{l=0}^{\infty} \left(
\sum_{m=0}^{l} \dfrac{[m]_q!}{m!} S(m,k) \binom lm_q
\beta_{l-m}^{(k)} ((1-x)_q, q) \right ) \dfrac{t^l}{[l]_q!}.
\endaligned
\eqno(30)
$$
Therefore, by (6) and (30), we obtain the following theorem,

\bigskip
{ \bf Theorem 10.}  For $k, l \in \mathbb{Z}_+$  and $ x \in [0,
1]$,  we have
$$ B_{k, l}(x, q) =
\dfrac{k! }{ [k]_q!} x^k \sum_{m=0}^{l} \dfrac{[m]_q!}{m!} S(m,k)
\beta_{l-m}^{(k)} ((1-x)_q, q)\binom lm_q ,$$ where
$\beta_{l}^{(k)} ((1-x)_q, q)$ are called the $l$-th $q$-Bernoulli
polynomials.
 \bigskip

From (15) and Theorem 10, we have the following corollary.

\bigskip
{ \bf Corollary 11.}  For $k, l \in \mathbb{Z}_+$  and $ x \in [0,
1]$, we have
$$ B_{k, l}(x, q) =
\dfrac{x^k }{ [k]_q!} \sum_{m=0}^{l} \dfrac{[m]_q!}{m!}\binom lm_q
\beta_{l-m}^{(k)} ((1-x)_q, q)\Delta^k 0^m .$$
 \bigskip

It is well known that
$$x^n= \sum_{k=0}^n \binom xk k! S(n, k),  \text{ (see [7])}. \eqno(31)$$

By (31) and Theorem 9, we easily see that
$$\sum_{k=0}^i \binom xk k! S(i, k)= \sum_{k=i}^n \dfrac{{\binom{k}{i}}_q }{{\binom{n}{i}}_q } B_{k, n}(x, q). $$

\bigskip

\begin{center} {\bf 3.  A matrix representation for $q$-Bernstein
polynomials}
\end{center}

Given a polynomial is written as a linear combination of
$q$-Bernstein basis functions:

$$B_q(x)= C_0^q B_{0, n}(x, q) + C_1^q B_{1, n}(x, q)+ \cdots + C_n^q B_{n, n}(x, q). \eqno(32)$$
It is easy to write (32) as a dot product of two vectors:

$$ B_q(x)= \left( \begin{array}{c} B_{0, n}(x, q),  B_{1, n}(x, q), \ldots,   B_{n, n}(x,
q)\end{array} \right)  \left( \begin{array}{c} C_0^q \\
C_1^q \\  \vdots \\ C_n^q  \end{array} \right). \eqno(33)$$

Now, we can convert (33) to
$$ B_q(x)= \left( \begin{array}{c} 1,  x, \ldots , x^n \end{array} \right)
 \begin{aligned} \begin{pmatrix}
b_{0,0}^q & 0 & \hdots  & 0    \\
b_{1,0}^q & b_{1,1}^q & \hdots  & 0   \\
\vdots & \vdots & \ddots  & \vdots   \\
b_{n,0}^q & b_{n,1}^q &   \hdots  & b_{n,n}^q
 \end{pmatrix} \end{aligned}
 \left( \begin{array}{c} C_0^q \\
C_1^q \\  \vdots \\ C_n^q  \end{array} \right),$$ where
$b_{i,j}^q$ are the coefficients of the power basis that are used
to determine the respective $q$-Bernstein polynomials.

From (5) and (6),  we note that

$$ \aligned  & B_{0,2}(x, q)=(1-x)_q^2 = \sum_{l=0}^2 \binom 2l_q (-1)^l  q^{ \binom l2} =1-[2]_qx +qx^2 \\
& B_{1,2}(x, q)= \binom 21_q x(1-x)_q =[2]_q x(1-x)=[2]_q x-[2]_q x^2 \\
&B_{2,2}(x, q)= x^2.
\endaligned
$$
In the quadratic case($n=2$), the matrix can be represented by
$$ B_q(x)= \left( \begin{array}{c} 1,  x, x^2 \end{array} \right)
 \begin{aligned} \begin{pmatrix}
1& 0 &  0    \\
-[2]_q & [2]_q &  0   \\
q & -[2]_q & 1
 \end{pmatrix} \end{aligned}
 \left( \begin{array}{c} C_0^q \\
C_1^q \\  C_2^q  \end{array} \right).$$

ACKNOWLEDGEMENTS. The present research has been conducted by the Research Grant of
Kwangwoon University in 2010.

\bigskip
\bigskip
\begin{center}{\bf REFERENCES}\end{center}
\begin{enumerate}

\item {M. Acikgoz, S. Araci,} {  A study on the integral of the product of several type
Bernstein polynomials,} { IST Transaction of Applied
Mathematics-Modelling and Simulation}, 2010.

\item {S. Bernstein,} {  D\.{e}monstration du th\.{e}or\.{e}me de Weierstrass, fonde\.{e} sur le Calcul des probabilities,}
 { Commun. Soc. Math, Kharkow (2)}, 13 (1912-1913), 1-2.

\item
{ I. N. Cangul, V. Kurt, H. Ozden, Y. Simsek}, { On the
higher-order $w$-$q$-Genocchi numbers,} { Adv. Stud. Contemp.
Math.}, {19} (2009), 39-57.

\item
{N. K. Govil, V. Gupta,} {Convergence of
$q$-Meyer-K\"{o}nig-Zeller-Durrmeyer operators,} { Adv. Stud.
Contemp. Math.}, {19} (2009), 97-108.

\item
{V. Gupta, T. Kim, J, Choi, Y.-H. Kim,} {Generating function for
$q$-Bernstein, $q$-Meyer-K\"{o}nig-Zeller and $q$-Beta basis,} {
Automation Computers Applied Mathematics}, {19} (2010), 7-11.

\item
{T. Kim,} { $q$-extension of the Euler formulae and trigonometric
functions,} { Russ. J. Math. Phys.},  { 14} (2007), 275-278.

\item
{ T. Kim,}  {   $q$-Volkenborn integration}, {Russ. J.  Math.
Phys.},  {  9} (2002), 288-299.

\item
{ T. Kim,}  {   $q$-Bernoulli numbers and polynomials associated
with Gaussian binomial coefficients},{ Russ. J. Math. Phys.}, {
15} (2008), 51-57.

\item
{T.Kim, L.C. Jang, H. Yi},{ Note on the modified $q$-Bernstein
polynomials}, Discrete Dynamics in Nature and Society(in press),
2010.

\item
{T. Kim, J. Choi, Y.-H. Kim,}  {  Some identities on the
$q$-Bernstein polynomials,  $q$-Stirling numbers and $q$-Bernoulli
numbers,} { Adv. Stud. Contemp. Math.}, {20} (2010), 335-341.

\item {T. Kim,}  {  Note on the Euler $q$-zeta functions,} {J.
Number Theory}, { 129} (2009), 1798-1804.

\item
{T. Kim,} { Barnes type multiple $q$-zeta function and $q$-Euler
polynomials,} { J. Phys. A: Math. Theor.}, {43} (2010) 255201,
11pp.

\item
{  V. Kurt}, { A further symmetric relation on the analogue of the
Apostol-Bernoulli and the analogue of the Apostol-Genocchi
polynomials,} { Appl. Math. Sci.}, 3 (2009), 53-56.

\item
{  B. A. Kupershmidt}, { Reflection symmetries  of $q$-Bernoulli
polynomials,} { J. Nonlinear  Math. Phys.}, 12 (2005), 412-422.

\item
{  G. M. Phillips}, { Bernstein polynomials based on the
$q$-integers,} { Annals of Numerical Analysis}, 4 (1997), 511-514.

\item
{  G. M. Phillips}, {On generalized  Bernstein polynomials,} {
Griffiths, D. F., Watson, G. A.(eds): Numerical Analysis},
Singapore: World Scientific, 263-269, 1996.

\item
{  Y. Simsek, M. Acikgoz}, { A new generating function of
$q$-Bernstein-type polynomials and their interpolation function,}
{ Abstract and Applied Analysis}, Article ID 769095 (2010), 12 pp.

 \item
{  L. C. Alberto}, { Probability and Random Processes for Electrical Engineering,}
{ Addison Wesley Longman}, 1994.

 \item
{  L. C. Biedenharn}, { The quantum group $SU_q(2)$ and a $q$-analogue of the boson operator,}
{ J. Phys. A}, 22(1989), L873-L878.

\item
{  S. C. Jing}, { The $q$-deformed binomial diftribution and its asymptotic behaviour,}
{ J. Phys. A}, 17(1994), 493-499.

\end{enumerate}

\end{document}